\input amstex
\magnification 1200
\documentstyle{amsppt}
\topmatter
\title
On degree of polynomial mappings of $\Bbb{C}^2$ to $\Bbb{C}^2$.
\endtitle
\author P. Katsylo
\endauthor
\address
Independent University of Moscow
\endaddress
\email katsylo\@ium.ips.ras.ru
\endemail
\date April 2, 1996
\enddate
\abstract
We prove two results about degree of polynomial
mappings of $\Bbb{C}^2$ to $\Bbb{C}^2$.
\endabstract
\endtopmatter
\document
\define\C{\Bbb C}
\define\A{\Bbb A}
\define\N{\Bbb N}
\define\Vn{{\widetilde V}_n}
\define\M{\widetilde{M}}
\define\LL{\widetilde{L}}
\bf \S0. \rm
\rm
Let $X_1, X_2$ be the canonical coordinates in $\C^2$,
$\C \lbrack X_1, X_2 \rbrack_{\le m}$ be
the linear space of polynomials of degree $\le m$
in the variable $X = (X_1, X_2)$.
Fix $n = (n_1, n_2) \in \N^2$ and put
$$
\Vn = \C [X_1, X_2]_{\le n_1} \times \C [X_1, X_2]_{\le n_2}.
$$
For $F = (F_1, F_2) \in \Vn$
consider the polynomial mapping
(denote it by the same letter)
$$
\gather
F : \C^2 \longrightarrow \C \times \C , \\
X = (X_1, X_2) \mapsto (F_1 (X_1, X_2), F_2 (X_1, X_2)) = F(X)
\endgather
$$
and let $J(F)$ be jacobian of the mapping $F$.
In the linear space $\Vn$
consider the following closed subvarieties
$$
\split
D_i = \{ & F \in \Vn \mid deg J(F) \le i \},
             \qquad 0 \le i \le n_1 + n_2 - 2, \\
W_i = \{ & F \in \Vn | dim F^{-1}(0) > 0 \ or \
             |F^{-1}(0)| \le i \ or \\
         & deg F_1 < n_1 \ and \ deg F_2 < n_2 \}, \qquad 0 \le i \le n_1
         n_2.  \endsplit $$
The results of the article are around the following problem.
\proclaim{Problem 0.1} For that $k, l$ the inclusion
$D_k \subset W_l$
holds?
\endproclaim
Note that the inclusion $D_0 \subset W_1$
for all $n$ is equivalent to the 2-dimensional
Jacobian Conjecture \cite{3}.
\proclaim{Theorem 0.2} For all $k \ge 0$ the inclusion
$$
D_k \subset W_{min\{n_1,n_2\}(k+1)}
$$
holds.
In particular if jacobian of a polynomial mapping
$F$ is identically equal to 1, then degree of $F$
does not exceed $min\{n_1, n_2\}$.
\endproclaim
By $C [X_1, X_2]_m$ denote
the linear space of homogeneous polynomials
of degree $m$ in the variable $X = (X_1, X_2)$.
Paar of elements $F \in \Vn , H' \in \C [X_1, X_2]_1$
are called general iff for $i = 1$
or for $i = 2$ the restriction of the polynomial $F_i(X)$
on the line $H'(X) = 0$
is a polynomial of degree $n_i$.
\proclaim{Theorem 0.3}
Let $F \in \Vn ,\ H' \in \C [X_1, X_2]_1$
be general elements, $|F^{-1}(0)| < \infty$.
Define the linear subspaces
$$
K = K(F,H'), K_i = K_i (F,H') \subset \C [X_1,X_2]_{\le n_1 + n_2 -1}
$$
in the following way
$$
\gather
K = \{F_1 Q_2 + F_2 Q_1 | Q_i \in \C [X_1, X_2]_{n_i - 1} \}, \\
K_0 = 0, \\
K_{i + 1} = (K + H' K_i) \cap \C [X_1, X_2]_{\le n_1 + n_2 - 2}.
\endgather
$$
Then
$$
\gather
K_0 \subset \dots \subset K_i \subset \dots
\subset \C[X_1, X_2]_{n_1 + n_2 - 2}, \\
2 dim K_i \ge dim K_{i-1} + dim K_{i+1}, \quad i \ge 1, \\
|F^{-1}(0)| = n_1 n_2 - dim K_\infty.
\endgather
$$
\endproclaim
We prove Theorem 0.2 in \S 1 and Theorem 0.3 in \S 2 - \S 4.

\bf \S1. \rm
In this section we prove Theorem 0.2.
\par
Let us remember some facts on Puiseux serieses.
\par
A Puiseux series (at the infinity)
is a convergent for $|t| > R$ series
$$
\alpha (t) = \sum_{i \le i_0} a_i (t^{\frac1d})^i,
$$
where $R > 0$.
Degree $deg_t \alpha (t)$ of the series $\alpha (t)$
is the greatest $\frac{i}{d}$ such that $a_i \neq 0$.
We write $\alpha \sim  \beta$ iff
$$
\beta (t) = \sum_{i \le i_0} a_i \theta^i (t^{\frac1d})^i,
$$
where $\theta \in \C$ is a $d$-th root of $1$.
The series $\alpha (t)$ is called reduced iff g.c.d.
of $\{i | a_i \neq 0 \}$ is equal to $1$.
Suppose $\alpha (t)$ is a reduced Puiseux series.
The series $\alpha (t)$ defines
$d$-valued analytical function
$$
\gather
\widetilde\alpha : \{t \in \C | |t| > R \}
\longrightarrow \C, \\
t \mapsto \alpha (t).
\endgather
$$
The number $d$ is called the denominator of $\alpha (t)$
(denote it by $den(\alpha)$).
\par
By $\C \{ t \}$ denote
the ring of Puiseux serieses.
We have the canonical differential operator
$$
\frac{d}{dt} : \C \{ t \} \longrightarrow \C \{ t \}.
$$
By $\C [X_2] \{ X_1 \}$ denote
the ring of polynomials in $X_2$
that coefficients are Puiseux serieses in $X_1$.
We have the canonical differential operators
$$
\frac{\partial}{\partial X_i} : \C[X_2] \{ X_1 \}
\longrightarrow \C[X_2] \{ X_1 \}, \quad i = 1,2.
$$
\par
An element $G \in \C [X_1, X_2]$ is called
proper with respect to $X_2$ iff
$$
G(X_1, X_2) = G_0 X_2^p + G_1 (X_1) X_2^{p-1} + \dots + G_p(X_1),
$$
where $0 \neq G_0 \in \C$.
\par
Let $G$ be proper with respect to $X_2$
polynomial in the variables $X_1, X_2$.
One can decompose the polynomial $G$ in
the ring $\C [X_2] \{ X_1 \}$.
Namely there exist Puiseux serieses
$\alpha_1 , \dots , \alpha_m$
such that
\roster
\item $den(\alpha_1) + \dots + den(\alpha_m) = deg_{X_2}G$,
\item
$$
G(X_1, X_2) = G_0 \prod_{1 \le l \le m, \alpha_l \sim \alpha}
(X_2 - \alpha (X_1)).  $$
\endroster
The Puiseux serieses
$\alpha_1, \dots , \alpha_m$
are called roots of $G$ with respect to $X_2$.
\par
Let $G$ be proper with respect to $X_2$
polynomial in the variables $X_1, X_2$.
There is well known procedure to construct
roots of $G$ with respect to $X_2$
(by means of Newton's polygons).
If we replace the polynomial $G$ by $G + c$, where
$c \in \C$, then the roots
(and sometimes the number of the roots) are changed.
The following fact is a corolary
of the procedure of the construction of roots.
\proclaim{Lemma 1.1} There exists an nonempty open in $\C$
subset $C = C(G)$ such that for $c \in C$,
roots $\alpha_1, \dots , \alpha_m$ of $G + c$
with respect to $X_2$, and
$$
G(X_1, X_2 + \alpha_l(X_1)) + c =
X_2 G_{l1}(X_1) + X_2^2 G_{l2}(X_1) + \dots ,
\quad 1 \le l \le m
$$
we have
$$
deg_{X_1}G_{l1}(X_1) \ge 0, \quad 1 \le l \le m.
$$
\endproclaim
\par
Let $G_1, G_2$ be polynomials in $X_1, X_2$
such that $G_1$ is proper with respect to $X_2$
and system of equations
$$
\left\{
  \aligned G_1(X_1, X_2) = 0 \\
           G_2(X_1, X_2) = 0
  \endaligned
\right.
\tag {1.1}
$$
has finitely many solutions.
Suppose $\alpha_1, \dots , \alpha_m$ is roots of the
polynomial $G_1$
with respect to $X_2$;
then the number of solutions (with multiplicities)
of the system (1.1) is equal to
$$
\sum_{1 \le l \le m} den(\alpha_l) deg_{X_1} G_2 (X_1, \alpha_l (X_1))
$$
(the Zeuthen formula \cite{4}).
\demo{Proof of Theorem 0.2}
Suppose $F = (F_1, F_2) \in D_k$ and let
$J = J(F)$ be jacobian of $F$.
We may assume that
$deg F_2 = n_2 \ge deg F_1 = n_1 = deg_{X_2}F_1$.
We have to prove that
$$
deg F \le n_1 (k + 1) \quad or \quad dim F(\C^2) < 2.
$$
Let $\alpha_1, \dots , \alpha_m$ be roots of
$F_1$ with respect to $X_2$.
It follows from the procedure of construction
of roots that
$$
deg_{X_1} \alpha_l(X_1) \le 1, \quad 1 \le l \le m. \tag{1.2}
$$
Set
$$
\gather
F_{1l}(X_1, X_2) = F_1(X_1, X_2 + \alpha_l(X_1))
       \in \C [X_2] \{ X_1 \}, \\
F_{1l}(X_1, X_2) = X_2 F_{1l1}(X_1) + X_2^2 F_{1l2}(X_1) + \dots ,
\quad 1 \le l \le m.
\endgather
$$
\par
Suppose $dim F(\C^2) = 2$.
We may assume that
\roster
\item
the number of the solutions (with multiplicities)
of the system of equations
$$
\left\{ \aligned
    F_1(X_1, X_2) = 0 \\
    F_2(X_1, X_2) = 0
        \endaligned
\right. \tag{1.3}
$$
is equal to degree of the mapping $F$,
\item
$deg_{X_1}F_{1l1}(X_1) \ge 0, \ 1 \le l \le m$ (see Lemma 1.1).
\endroster
\par
Fix $1 \le l \le m$ and consider
$$
\gather
F_{2l}(X_1, X_2) = F_2 (X_1, X_2 + \alpha_l (X_1))
      \in \C [X_2] \{ X_1 \}, \\
F_{2l}(X_1, X_2) = F_{2l0} (X_1) + X_2 F_{2l1} (X_1) + \dots , \\
J_l (X_1, X_2) = J(X_1, X_2 + \alpha_l (X_1))
      \in \C [X_2] \{ X_1 \}, \\
J_l (X_1, X_2) = J_{l0}(X_1) + X_2 J_{l1} (X_1) + \dots .
\endgather
$$
From (1.2) it follows that
$$
deg_{X_1}J_{l0}(X_1) \le k. \tag{1.4}
$$
We have
$$
det\left(\left(\frac{\partial F_{il}}{\partial X_j}\right)_
{1 \le i,j \le 2}\right) (X_1, X_2) =
J_l (X_1, X_2)
$$
whence
$$
deg_{X_1}\left(\frac{d}{dX_1} F_{2l0} (X_1)\right) +
deg_{X_1} F_{1l1} (X_1) =
deg_{X_1} J_{l0} (X_1).
$$
Using (2) and (1.4), we get
$$
deg_{X_1} F_{2l0} (X_1) \le deg_{X_1} J_{l0} (X_1) + 1 \le k + 1. \tag{1.5}
$$
\par
The number of the solutions (with multiplicities)
of the system of equations (1.3) is equal to
$$
\gather
\sum_{1 \le l \le m} den(\alpha_l) deg_{X_1} F_2 (X_1, \alpha_l(X_1)) =
    \sum_{1 \le l \le m} den(\alpha_l) deg_{X_1} F_{2l0} (X_1) \le \\
\le \sum_{1 \le l \le m} den(\alpha_l) (k + 1) = n_1 (k + 1)
\endgather
$$
(see (1.5)). Therefore, degree of the mapping $F$
does not exceed
$min \{ n_1, n_2 \}(k+1)$.
\enddemo
\par
\bf \S 2. \rm
We use some standard facts of representation theory \cite{1}.
\par
For a linear space $V$ and dual space $V^\ast$ by
$$
\gather
\langle \ .\  ,\ .\ \rangle : V \times V^\ast
\longrightarrow \C, \\
(v, v^\ast) \mapsto \langle v, v^\ast \rangle
\endgather
$$
denote the canonical bilinear mapping.
\par
The group $SL_3$ acts canonically in the spaces
$\C^3, \C^{3\ast}$, $S^m\C^{3\ast}, \dots $.
Let $e_1, e_2, e_3$ be the standard basis of $\C^3$,
$x_1, x_2, x_3$ be the standard basis of $\C^{3\ast}$.
Set
$$
\Delta = \sum \frac{\partial}{\partial e_i} \otimes
              \frac{\partial}{\partial x_i} :
\C [\C^{3\ast}] \otimes \C[\C^3] \longrightarrow
\C [\C^{3\ast}] \otimes \C[\C^3].
$$
Recall that
$$
dim S^m\C^{3\ast} = \frac12 (m+1)(m+2).
$$
For $g \in S^m\C^{3\ast}$ put
$$
V(g) = \{ \overline{a} \in P\C^3 | g(\overline{a}) = 0 \}
\subset P\C^3.
$$
\par
Let $n_1, n_2$ be natural numbers.
Set  $n = (n_1, n_2)$, $N = n_1 n_2$, and
$$
V_n = S^{n_1}\C^{3\ast} \times S^{n_2}\C^{3\ast}.
$$
In this section we make the following.
\roster
\item
We define the covariant
$$
Q : V_n \longrightarrow \C
$$
(eliminant) and prove some of its properties.
\item
For $h \in \C^{3\ast}, h \neq 0, f \in V_n$ we define
affine space $\A (h) \subset P\C^3$
and the polynomial mapping
$$
I_h(f) : \A (h) \longrightarrow \C \times \C.
$$
We have $Q(f) = 0$ iff
$dim (I_h(f)^{-1} (0)) > 0$
or polynomial degree of $I_h(f)$ is less $(n_1; n_2)$.
\item
For $f \in V_n, h, h' \in \C^{3\ast}$ such that
$Q(f) \neq 0$, $h \neq 0$, $h' \neq 0$,
and $V(h) \cap V(h') \cap V(f_1) \cap V(f_2) = \emptyset$
we prove that
$$
|I_h(f)^{-1}(0)| = deg_tR(f, h' + th),
$$
where
$$
R : V_n \times \C^{3\ast} \longrightarrow \C
$$
is the resultant.
\endroster
\par
1.
Consider the resultant
$$
R : V_n \times \C^{3\ast} \longrightarrow \C.
$$
The resultant $R$ is a polyhomogeneous
(of polydegree $(n_1, n_2, n_1 n_2)$) invariant.
The resultant $R$ defines canonically
the polyhomogeneous (of polydegree $(n_2, n_1)$) covariant
$$
Q : V_n \longrightarrow (S^N \C^{3\ast})^\ast = S^N \C^3.
$$
We have
$$
R(f,h) = \langle Q(f), h^N \rangle
$$
for $(f,h) \in V_n \times \C^{3\ast}$.
It follows from this formula that
\roster
\item
if $f \in V_n$, then $Q(f) = 0$ iff
$dim(V(f_1) \cap V(f_2)) > 0 ,$
\item
if $f \in V_n,\ dim(V(f_1) \cap V(f_2)) = 0$, then
$Q(f) = l_1 \cdot \dots \cdot l_N,$
where $l_1, \dots , l_N \in \C^3$ and
$V(f_1) \cap V(f_2) = \{ \overline{l_1}, \dots , \overline{l_N} \}.$
\endroster
\par
2.
Supose $h \in \C^{3\ast},\ h \neq 0$.
Define the 2-dimensional affine space
$$
\A (h) = P\C^3 \setminus V(h).
$$
Let $\C[\A(h)]_{\le m}$ be the linear space
of polynomial mappings (of polynomial degree $\le m$)
of the affine space $\A(h)$ to $\C$.
The linear space
$\C[\A(h)]_{\le n_1} \times \C[\A(h)]_{\le n_2}$
is the space of polynomial mappings
(of polynomial degree $\le (n_1; n_2)$)
of the affine space $\A(h)$ to $\C \times \C$.
Fix the isomorphisms of the liner spaces
$$
\gather
i_h : S^m \C^{3\ast} \longrightarrow \C[\A(h)]_{\le m}, \\
i_h(g)(a) = \frac{g(a)}{h(a)^m}, \\
I_h : V_n \longrightarrow
\C[\A(h)]_{\le n_1} \times \C[\A(h)]_{\le n_2}, \\
I_h(f)(a) =
\left( \frac{f_1(a)}{h(a)^{n_1}}, \frac{f_2(a)}{h(a)^{n_2}} \right).
\endgather
$$
From (1) of item 1 it follows that
$Q(f) = 0$ iff
$dim(I_h(f)^{-1}(0)) > 0$ or polynomial degree
of $I_h(f)$ is less $(n_1; n_2)$.
\par
3.
\proclaim{Lemma 2.1} Suppose $f \in V_n$, $h, h' \in \C^{3\ast}$,
$Q(f) \neq 0$, $h \neq 0$, $h' \neq 0$, and
$V(h) \cap V(h') \cap V(f_1) \cap V(f_2) = \emptyset$;
then
$$
|I_h(f)^{-1}(0)| = deg_tR(f, h' + th).
$$
\endproclaim
\demo{Proof}
We have
$$
Q(f) = l_1 \cdot \dots \cdot l_N,
$$
where $l_i \in \C^3$. We may assume that
$h(l_i) \neq 0$ for $1 \le i \le d$
and $h(l_i) = 0$ for $d+1 \le i \le N$.
From the suppositions of the Lemma it follows that
$h'(l_i) \neq 0$ for $d+1 \le i \le N.$
We have
$$
\gather
I_h(f)^{-1}(0) = \{ \overline{l_1}, \dots , \overline{l_d} \}, \\
R(f, h' + th) = \langle Q(f), (h' + th)^N \rangle = \\
= \langle l_1 \cdot \dots \cdot l_N,
{h'}^N + tN{h'}^{(N-1)}h + \dots + t^Nh^N \rangle, \\
\langle l_1 \cdot \dots \cdot l_N, {h'}^{(N-d)}h^d \rangle \neq 0, \\
\langle l_1 \cdot \dots \cdot l_N, {h'}^{(N-i)}h^i \rangle = 0, \quad
d+1 \le i \le N.
\endgather
$$
From these formulas it follows the Lemma.
\enddemo
\par
\bf \S 3. \rm
We use the notations of \S 2.
\par
Suppose $n = (n_1, n_2) \in \N^2$.
\par
Recall a procedure of a calculation of the resultant
$$
R : V_n \times \C^{3\ast} \longrightarrow \C.
$$
Fix $(f,s) \in V_n \times \C^{3\ast}$ and
consider the complex of vector spaces
$$
M(f,s): \quad 0 \longrightarrow M @> \beta(f,s)>>
M' @> \beta'(f,s) >> M'' \longrightarrow 0,
$$
where
$$
\gather
M = S^{n_1 - 2}\C^{3\ast} \times S^{n_2 - 2}\C^{3\ast}, \\
M' = S^{n_1 - 1}\C^{3\ast} \times S^{n_2 - 1}\C^{3\ast}
     \times S^{n_1 + n_2 - 2}\C^{3\ast}, \\
M'' = S^{n_1 + n_2 - 1}\C^{3\ast}, \\
\beta(f,s)(r_1, r_2) = (s r_1, s r_2, f_1 r_2 + f_2 r_1), \\
\beta'(f,s)(q_1, q_2, g) = (f_1 q_2 + f_2 q_1 - s g).
\endgather
$$
The determinant of the complex $M(f,s)$ is
equal to $R(f,s)$.
\par
One can calculate the determinant of the complex
$M(f,s)$ in the following way \cite{2}.
Fix $a \in \C^3, s(a) \neq 0$
and consider the linear mappings
$$
\gather
\alpha(a) : M' \longrightarrow M, \\
(q_1, q_2, g) \mapsto (\Delta (q_1 a), \Delta (q_2 a)), \\
t \alpha(a) \circ \beta(f,s) : M \longrightarrow M, \\
(t \alpha(a), \beta'(f,s)) : M' \longrightarrow M \times M''.
\endgather
$$
We have
$$
R(f,s) = det(t \alpha(a), \beta'(f,s))
(det(t \alpha(a) \circ \beta(f,s)))^{-1}
$$
(we calculate $det$ with respect to the canonical
$SL_3$-invariant forms of maximal degree in the
spaces $S^m\C^{3\ast}$).
\proclaim{Lemma 3.1}
$$
det(t \alpha(a) \circ \beta(f,s)) = c(ts(a))^{dimM},
$$
where $0 \neq c \in \C$.
\endproclaim
\demo{Proof}
The function $det(\alpha(a) \circ \beta(f,s))$
is a nonzero polyhomogeneous (of polydegree $(dim M, dim M)$)
$SL_3$-invariant function in $(a,s) \in \C^3 \times \C^{3\ast}$.
As is known, there exists only one
(up to a nonzero factor) that function
that is $(s(a))^{dimM}$.
\enddemo
Therefore,
$$
R(f,s) = det(t \alpha(a), \beta'(f,s))
c^{-1}(ts(a))^{-dimM}. \tag{3.1}
$$
We need the following fact of linear algebra.
\proclaim{Lemma 3.2}
Let
$$
\eta, \eta' : V \longrightarrow V'
$$
be linear mappings, $dim V = dim V'$, and
$det(\eta' + t\eta) \neq 0$. Define
$$
L_i = L_i(\eta,\eta') \subset V', \quad i \ge 0
$$
in the following way
$$
\gather
L_0 = 0, \\
L_{i+1} = \eta'(\eta^{-1}(L_i)) \cap Im \eta.
\endgather
$$
Then
$$
\gather
L_0 \subset \dots \subset L_i \subset \dots \subset Im \eta, \\
2 dim L_i \ge dim L_{i-1} + dim L_{i+1}, \quad i \ge 1, \\
deg_tdet(\eta' + t\eta) = dim V - dim \eta^{-1}(0) - dim L_\infty.
\endgather
$$
\endproclaim
\demo{Proof}
Let $b_1, \dots , b_n$ be a basis of $V$,
$m \ge 0$, and
$$
\epsilon : V' \longrightarrow V
$$
be an isomorphism such that
$$
\gather
(\epsilon \circ \eta')(b_i) =
  \left\{\aligned &b_i \quad for \ 0 \le i \le m, \\
                &0   \quad for \ m+1 \le i \le n,
         \endaligned
  \right. \\
  (\epsilon \circ \eta)(\langle b_1, \dots , b_m \rangle )
  \subset \langle b_1, \dots , b_m \rangle , \\
  (\epsilon \circ \eta)(\langle b_{m+1}, \dots , b_n \rangle )
  \subset \langle b_{m+1}, \dots , b_n \rangle ,
\endgather
$$
the matrix of $\epsilon \circ \eta$
in the basis $b_1, \dots , b_n$
is a Jordan matrix.
It can easily be checked that
the Lemma holds for
$\epsilon \circ \eta$, $\epsilon \circ \eta'$.
Therefore, the Lemma holds for $\eta, \eta'$.
\enddemo
\proclaim{Theorem 3.3}
Suppose $f \in V_n$, $a \in \C^3$,
$h, h' \in \C^{3\ast}$,
$Q(f) \neq 0$, $h \neq 0$, $h' \neq 0$,
$(h' + t h)(a) \neq 0$, and
$V(h) \cap V(h') \cap V(f_1) \cap V(f_2) = \emptyset$ and
define
$$
L_i = L_i (f,a,h,h') \subset M \times M'', \quad i \ge 0
$$
in the following way
$$
\gather
L_0 = 0, \\
L_{i+1} = (0,\beta'(f,h'))((\alpha(a),\beta'(0,h))^{-1}(L_i))
          \cap Im(\alpha(a), \beta'(0,h));
\endgather
$$
then
$$
\gather
L_0 \subset \dots \subset L_i \subset \dots \subset
\{ 0 \} \times M'', \\
2 dim L_i \ge dim L_{i-1} + dim L_{i+1}, \quad i \ge 1, \\
deg_tR(f, h' + th) = n_1 n_2 - dim L_\infty.
\endgather
$$
\endproclaim
\demo{Proof} The Theorem is a corolary of
(3.1), Lemma 3.2,
the formula for dimension of $S^m\C^{3\ast}$ (see \S 2),
and the following evident fact
$$
dim(\alpha(a), \beta'(0,h))^{-1}(0) = n_1 + n_2.
$$
\enddemo
\par
\bf \S 4. \rm
In this section we prove Theorem 0.3.
We use the notations of \S 2 and \S3.
\par
Fix the isomorphisms of the linear spaces
(denote them by one letter)
$$
\gather
\theta : \C[X_1, X_2]_{\le m} \longrightarrow
         S^m\C^{3\ast} , \quad m \ge 0, \\
G(X_1, X_2) \mapsto x_3^m G(\frac{x_1}{x_3},\frac{x_2}{x_3}).
\endgather
$$
We see that
$$
\theta = i_{x_3}^{-1} \circ j_3^\ast,
$$
where
$$
\gather
j_3 : \A (x_3) \longrightarrow \C^2, \\
(a_1: a_2: 1) \mapsto (a_1, a_2)
\endgather
$$
is the isomophism of the affine varieties.
\par
Set
$$
\gather
\M = \C[X_1, X_2]_{\le n_1 - 2} \times \C[X_1, X_2]_{n_2 - 2}, \\
\M' = \C[X_1, X_2]_{\le n_1 - 1} \times \C[X_1, X_2]_{n_2 - 1}
      \times \C[x_1, X_2]_{n_1 + n_2 - 2}, \\
\M'' = \C[X_1, X_2]_{n_1 + n_2 - 1}.
\endgather
$$
The isomorphisms $\theta$
define canonically the following isomorphisms
(denote them by the same letter)
$$
\theta: \quad \Vn \rightarrow V_n, \quad \M \rightarrow M,
\quad \M' \rightarrow M', \quad \M'' \rightarrow M''.
$$
\demo{Proof of Theorem 0.3}
Set $f = \theta (F)$, $h' = \theta (H')$.
From the suppositions of the Theorem
it follows that $Q(f) \neq 0$,
$(h' + tx_3)(e_3) \equiv t$, and
$V(x_3) \cap V(h') \cap V(f_1) \cap V(f_2) = \emptyset$.
Using Lemma 2.1 and Theorem 3.3, we get
$$
\gather
|F^{-1}(0)| = |I_{x_3}(f)^{-1}(0)| = deg_tR(f,h' + tx_3) = \\
n_1n_2 - dim L_\infty (f, e_3, x_3, h').
\endgather
$$
\par
The following linear mappings correspond to
$(0,\beta'(f,h'))$ and $(\alpha(e_3), \beta'(0,x_3))$
under the isomorphisms $\theta$:
$$
\gather
\gamma'(F,H') = \theta^{-1} \circ (0,\beta'(f,h')) \circ \theta :
       \M' \longrightarrow \M \times \M'', \\
(Q_1, Q_2, G) \mapsto (0, F_1Q_2 + F_2Q_1 - H'G), \\
\gamma(F,H') =
       \theta^{-1} \circ (\alpha(e_3), \beta'(0,x_3)) \circ \theta :
       \M' \longrightarrow \M \times \M'', \\
(Q_1, Q_2, G) \mapsto
       (\Delta_{n_1 -1}(Q_1),\Delta_{n_2 -1}(Q_2), -G),
\endgather
$$
where
$$
\gather
\Delta_m : \C[X_1, X_2] \longrightarrow \C[X_1, X_2], \\
X_1^{m_1} X_2^{m_2} \mapsto (m-m_1-m_2) X_1^{m_1} X_2^{m_2}.
\endgather
$$
\par
Define the linear subspaces
$$
\LL_i = \LL_i(F,H') \subset \M \times \M''
$$
in the following way
$$
\gather
\LL_0 = 0, \\
\LL_{i+1} = \gamma'((F,H')(\gamma(F,H')^{-1}(\LL_i))
        \cap Im \gamma(F,H').
\endgather
$$
We have
$$
\LL_i (F, H') = \theta (L_i (f, e_3, x_3, h')), \quad i \ge 0
$$
and therefore,
$$
\gather
\LL_0 \subset \dots \subset \LL_i \subset \dots
      \subset \M \times \M'', \\
2 dim \LL_i \ge dim \LL_{i-1} + \LL_{i-1},
\quad i \ge 1, \tag {4.1} \\
|F^{-1}(0)| = n_1 n_2 - dim \LL_\infty.
\endgather
$$
It follows from the definitions that
$$
\gather
\LL_i \subset \{ 0 \} \times \M'', \qquad \gamma(F,H')^{-1}(0,Z) = \\
= \left\{ \aligned
        &0 \quad if \ degZ = n_1 + n_2 -1, \\
        &\C[X_1,X_2]_{n_1 - 1} \times \C[X_1,X_2]_{n_2 - 1}
        \times \{ -Z \} \quad if \ degZ \le n_1 + n_2 - 2.
     \endaligned \right.
\endgather
$$
It is easy to prove by induction that
$$
\LL_i (F,H') = \{ 0 \} \times K_i(F,H')
$$
and therefore,
$$
\gather
K_0 \subset \dots \subset K_i \subset \dots \subset
          \C[X_1,X_2]_{n_1 + n_2 - 2}, \\
2 dim K_i \ge dim K_{i-1} + dim K_{i+1}, \quad i \ge 0, \\
dim \LL_\infty = dim K_\infty, \\
|F^{-1}(0)| = n_1 n_2 - dim K_\infty.
\endgather
$$
(see (4.1)).
\enddemo


\Refs
\ref \no1
\by W.Fulton and J.Harris
\book Representation Theory
\publ Springer Graduate Text in Math., v. 129,
      Springer-Verlag, Berlin
\yr 1991
\endref
\ref \no2
\by I.M.Gelfand, M.M.Kapranov, and A.V.Zelevinsky
\book Discriminants, Resultants and Multidimensional Determinants
\publ Birkauser, Boston
\yr 1994
\endref
\ref \no3
\by H.Bass, E.H.Connell, and D.Wright
\paper The Jacobian Conjecture: reduction of degree and formal
       expansion of the inverse
\jour Bull. Amer. Math. Soc.
\vol 1982 \issue 7 \pages 287 - 330
\endref
\ref \no4
\by W.Fulton
\book Intersection Theory
\publ Springer-Verlag Berlin, Heidelberg
\yr 1984
\endref
\endRefs
\enddocument